\documentclass[12pt,letterpaper]{article}

\usepackage{fullpage}
\usepackage{amssymb}
\usepackage{amsmath}
\usepackage{textcomp}
\usepackage{float}

\def\cent{\textcent}
\def \mod#1 #2{#1\ ({\rm mod}\ #2)}
\def \nodiv{{|\kern-3.9pt/}}

\def \prend{\vrule depth-1pt height7pt width6pt}
\def \proof{\bigbreak\noindent{\bf Proof.\ \ }}
\def \endpf{{\ \ \prend \medbreak}}

\def \Reals{\mathbb R}

\def \Zee{\mathbb Z}

\title{The Computational Complexity of the Local Postage Stamp Problem}

\author{Jeffrey Shallit\thanks{Currently on sabbatical (until July 2002)
at Department of Computer Science, University of Arizona,
P.~O. Box 210077, Tucson, AZ  85721-0077, USA.}\\
Department of Computer Science\\
University of Waterloo\\
Waterloo, Ontario, Canada  N2L 3G1\\
{\tt shallit@graceland.uwaterloo.ca} 
}

\newtheorem{theorem}{Theorem}
\newtheorem{lemma}[theorem]{Lemma}

\begin{document}

\maketitle

%

\begin{abstract}
The well-studied local postage stamp problem (LPSP) is the following:  given
a positive integer $k$,
a set of positive integers $1= a_1 < a_2 < \cdots < a_k$
and an integer $h \geq 1$, what is the smallest positive integer which
{\it cannot\/} be represented as a linear combination
$\sum_{1 \leq i \leq k} x_i a_i$ where $\sum_{1 \leq i \leq k} x_i \leq h$
and each $x_i$ is a non-negative integer?
In this note we prove that LPSP is NP-hard under Turing reductions,
but can be solved in polynomial time if $k$ is fixed.
\end{abstract}

\section{Introduction}

The {\sl local postage-stamp problem}, or LPSP for short, can be
informally defined as follows.  One is given a supply of stamps
of $k$ different denominations, $1 = a_1 < a_2 < \cdots < a_k$,
and an envelope that has room for at most $h$ different stamps.  What
is the smallest amount of postage $N_h( a_1, a_2, \ldots, a_k)$
that {\it cannot\/} fit on the envelope?
For example, if the available denominations are 1\cent, 4\cent,
7\cent, and 8\cent, and the envelope has room for $3$ stamps,
then all amounts of postage $\leq$ 24\cent\ can be provided 
but 25\cent\ cannot.  Hence $N_3 (1,4,7,8) = 25$.

A more formal statement of LPSP is given in the abstract.

In this note we consider the computational complexity of LPSP.
If $N = \sum_{1 \leq i \leq k} x_i a_i$, we call 
$(x_1, x_2, \ldots, x_k)$ a {\sl representation} for $N$ and
$\sum_{1 \leq i \leq k} x_i$ the {\sl weight} of the representation.
If further $\sum_{1 \leq i \leq k} x_i$
is minimum among all representations for $N$, we call $(x_1, x_2, \ldots, x_k)$ a 
{\sl minimum-weight representation} for $N$.
If the denominations $a_i$ and bound $h$ are given
in {\it unary}, then a simple dynamic programming algorithm can determine
the minimum-weight representation for all integers $N \leq ha_k +1$
in polynomial-time, and hence we can compute
$N_h(a_1, a_2, \ldots, a_k)$ in polynomial time.
We therefore assume for the rest of this paper that
all inputs are provided in {\it binary}.

LPSP was apparently introduced by Rohrbach \cite{Rohrbach:1937a,Rohrbach:1937b}
in 1937, and since then dozens of papers have been written about it
and a variant, the {\sl global postage-stamp problem} (GPSP);
see Guy \cite[pp.\ 123--127]{Guy:1994}) for a brief survey.  
Despite this, no general results on the computational complexity of the
problem seem to be known up to now;
for example, Alter and Barnett \cite{Alter&Barnett:1980} asked if 
$N_h (a_1, \ldots, a_k)$ ``can be expressed by a simple formula''.
Selmer \cite{Selmer:1980} discussed efficient algorithms for
the case where $k \leq 3$.  

In the next section we prove that LPSP is NP-hard under Turing reductions,
and in Section 3 we give a polynomial-time algorithm for LPSP when $k$ is fixed.

\section{LPSP is NP-hard}

We prove that LPSP is NP-hard by reducing from a related problem,
the {\sl Frobenius problem} (see,
for example, Guy \cite[pp.\ 113--114]{Guy:1994}).
In the Frobenius problem,
we are given an integer $k \geq 1$ and $k$ positive integers
$a_1, a_2, \ldots, a_k$ with $\gcd(a_1, a_2, \ldots, a_k) = 1$, and we are asked
to compute $g(a_1, a_2, \ldots, a_k)$,
the {\it largest\/} integer which {\it cannot\/}
be expressed as a non-negative integer linear combination
$\sum_{1 \leq i \leq k} x_i a_i$.
The Frobenius problem is well-studied, but it was only fairly
recently that it was proved NP-hard (under Turing reductions)
by Ram\'{\i}rez-Alfons\'{\i}n \cite{Ramirez-Alfonsin:1996}.

Before we give the reduction, we need a technical lemma.

\begin{lemma}
Let $1 = a_1 < a_2 < \cdots < a_k$.    
Define 
$$h_0 = \sum_{1 \leq i \leq k} \left\lfloor {{a_{i+1}} \over {a_i}} \right\rfloor,$$
and
$$h_1 = h_0 + \left\lceil  {{(h_0+1) a_{k-1}} \over {a_k - a_{k-1}}} \right\rceil.$$
Then
	\begin{itemize}
	\item[(a)] $N_{h_0} (a_1, a_2, \ldots, a_k) > a_k$;

	\item[(b)] $N_{h_0+i} (a_1, a_2, \ldots, a_k) > (i+1) a_k$ for all $i \geq 0$;

	\item[(c)] $N_h (a_1, a_2, \ldots, a_k) > (h+1) a_{k-1} - a_k$ for all
	$h \geq h_1$.

	\item[(d)] $N_{h+1} (a_1, a_2, \ldots, a_k) = N_h(a_1, a_2,\ldots, a_k) + a_k$
for all $h \geq h_1$.

	\item[(e)]  There exists a constant $c \geq -1$ such that
	$h a_k - N_h (a_1, a_2, \ldots, a_k) = c$ for all $h \geq h_1$.

	\item[(f)] If $h \geq h_1$, then 
$N_h (a_1, a_2, \ldots, a_k)
= h a_k - g(a_k - a_{k-1}, a_k - a_{k-2}, \ldots, a_k - a_1, a_k)$,
where $g$ is the Frobenius number.
	\end{itemize}
\label{selmer}
\end{lemma}

\noindent {\bf Remark.}  Parts (a)-(f) can be essentially found more or
less verbatim in the paper of Selmer
\cite{Selmer:1980}; the only difference in our presentation is that we explicitly
compute the bounds $h_0, h_1$.

\proof
     (a)  Consider finding a representation $ n = \sum_{1 \leq i \leq k} x_i a_i$
for an integer $n$, $0 \leq n \leq a_k$, using the greedy algorithm.  We use
at most $\lfloor a_k/a_{k-1} \rfloor$ copies of $a_{k-1}$, then at most
$\lfloor a_{k-1}/a_{k-2} \rfloor$ copies of $a_{k-2}$, etc.  The choice of $h_0$
thus allows us to form the greedy representation of all such $n$.

     (b)   We prove this by induction on $i$.  For $i = 0$ the result 
is just part (a).  Otherwise, suppose
$N_{h_0+i} (a_1, a_2, \ldots, a_k) > (i+1) a_k$.  Then every integer
$n$, $0 \leq n \leq (i+1)a_k$ has a representation of weight $\leq h_0 + i$ and
by adding a single copy of $a_k$, we can represent every integer
$m$, $0 \leq m \leq (i+2)a_k$ with weight $\leq h_0 + i + 1$.

     (c)  Set $i = j + \lceil {{(h_0+1) a_{k-1}} \over {a_k - a_{k-1}}} \rceil$,
$j \geq 0$, and apply (b).

     (d)  The numbers that have representations of weight $\leq h+1$ can be
divided into two not necessarily disjoint subsets:
\begin{eqnarray*}
S_1 &=& \lbrace n \ : \ \text{some representation of $n$ of weight $\leq h+1$
has $x_k > 0$} \rbrace \\
S_2 &=& \lbrace n \ : \ \text{some representation of $n$ of weight $\leq h+1$
has $x_k = 0$} \rbrace.
\end{eqnarray*}
Now every element of $S_1$ can be written as $a_k + t$, where $t$ has
a representation of weight $\leq h$.  It follows that 
$$\lbrace a_k, a_k + 1, a_k + 2,  \ldots, N_h(a_1,a_2,\ldots, a_k) +a_k -1 \rbrace 
\subseteq S_1,$$
but $N_h(a_1, h_2, \ldots, a_k) + a_k \not\in S_1$.  
On the other hand, the numbers in $S_2$ have representations of weight $\leq h+1$
using just the numbers $\lbrace a_1, a_2, \ldots, a_{k-1} \rbrace$, and
so the largest element of $S_2$ is $\leq (h+1) a_{k-1}$.  
Furthermore, by (a) and the fact that $h_1 \geq h_0$,
we have $\lbrace 0, 1, \ldots, a_k \rbrace \subseteq S_2$.
It follows that provided 
\begin{equation}
N_h(a_1, a_2, \ldots, a_k) + a_k > (h+1) a_{k-1},
\label{ineq}
\end{equation}
we have $N_{h+1} (a_1, a_2, \ldots, a_k) = N_h (a_1, a_2, \ldots, a_k) + a_k$.
But (\ref{ineq}) follows from (c).

     (e)  Using (d), a simple induction gives
$N_{h_1 + i} (a_1, a_2, \ldots, a_k) = N_{h_1} (a_1, a_2 , \ldots, a_k) + i a_k$
for all $i \geq 0$.
Then $(h_1 + i) a_k - N_{h_1 + i} (a_1, a_2, \ldots, a_k) 
= h_1 a_k - N_{h_1} (a_1, a_2, \ldots, a_k)$ for all $i \geq 0$;
so we may take $c = h_1 a_k - N{h_1} (a_1, a_2, \ldots, a_k)$.
Since $N_{h_1} (a_1, a_2, \ldots, a_k) \leq h_1 a_k + 1$, it follows that $c \geq -1$.

     (f)  Suppose $h \geq h_1$.  Then by (e)
we have $h a_k - N_h (a_1, a_2, \ldots, a_k)$
is independent of $h$.  Now $h a_k - t$ has a representation
of weight $h' \leq h$ iff
\begin{eqnarray*}
t &=& h a_k - \sum_{1 \leq i \leq k} x_i a_i \\
&=& \left( \sum_{1 \leq i \leq k-1} x_i (a_k - a_i) \right) + (h-h') a_k,
\end{eqnarray*}
i.e., if $t$ has a representation of {\it any\/} weight
using the basis
$a_k - a_{k-1}, a_k - a_{k-2}, \ldots, a_k - a_1, a_k$,
since we can choose $h$ to be arbitrarily large.
But the largest $t$ with no representation in the basis 
$a_k - a_{k-1}, a_k - a_{k-2}, \ldots, a_k - a_1, a_k$
is just the Frobenius number
$$g(a_k - a_{k-1}, a_k - a_{k-2}, \ldots, a_k - a_1, a_k).$$
(Since
$a_1 = 1$, we have $\gcd(a_k - a_{k-1}, a_k - a_{k-2}, \ldots, a_k - a_1, a_k ) = 1$.)
It follows that 
$h a_k - N_h (a_1, a_2, \ldots, a_k)
= g(a_k - a_{k-1}, a_k - a_{k-2}, \ldots, a_k - a_1, a_k)$.
\endpf

\begin{theorem}
     Given positive integers $b_1 < b_2 < \cdots < b_k$ with 
$\gcd(b_1, b_2, \ldots, b_k) = 1$, we can determine in polynomial
time integers $h, a_1 =1, a_2, \ldots, a_k, a_{k+1}, a_{k+2}$ such that
$g(b_1, b_2, \ldots, b_k) = h a_{k+2} - N_h(a_1, a_2, \ldots, a_k, a_{k+1}, a_{k+2})$.
\end{theorem}

\proof
By a theorem of Brauer 
\cite[Corollary to Thm. 1]{Brauer:1942}, 
we know that 
$g(b_1, b_2, \ldots, b_k) < b_k b_1$.
Define $b_{k+1} = b_k b_1$
and $b_{k+2} = b_k b_1 + 1$.  Then clearly
$g(b_1, b_2, \ldots, b_k) = g(b_1, b_2, \ldots, b_k, b_{k+1}, b_{k+2})$.
Now we have, by Lemma~\ref{selmer} (f), that
$$h_1 a_{k+2} - N_{h_1} (b_{k+2} - b_{k+1}, b_{k+2} - b_k, \ldots, b_{k+2} - b_1, b_{k+2})
= g(b_1, b_2, \ldots, b_k, b_{k+1}, b_{k+2}),$$
where
$$h_1 = h_0 + \left\lceil  {{(h_0+1) a_{k+1}} \over {a_{k+2} - a_{k+1}}} \right\rceil,$$
and
$$h_0 = \sum_{1 \leq i \leq k+2} \left\lfloor {{a_{i+1}} \over {a_i}} \right\rfloor.$$
Note that $b_{k+2} - b_{k+1} = 1$.
\endpf

    Since the Frobenius problem reduces to LPSP, and the Frobenius problem
is NP-hard, so is LPSP.

\section{A polynomial-time algorithm for fixed $k$}

I observe that results of 
Kannan \cite{Kannan:1989,Kannan:1990} provide
a polynomial-time algorithm for the local postage-stamp problem
for any fixed dimension $k$.  
Let $Q$ be a given copolyhedron (an intersection of a finite number
of half-spaces, possibly closed, possibly open) in $\Reals^{p+l}$.
Let $A$ be an $m \times n$ matrix, $B$ be a
$m \times p$ matrix, and $C$ be a column vector of dimension $m$,
all with integer entries.
Kannan proved that assertions of the form

\centerline{``$\forall y \in Q/\Zee^l \
\ \exists x \in \Zee^n \ \ \text{such that} \ Ax + By \leq C$''}

\noindent can be tested in polynomial time when $l+n+p$ is fixed.  
(Here $Q/\Zee^l = \lbrace y \in \Reals^p \ :
\ \text{there exists}\ t\in \Zee^l \ \text{such that}\ [y, t] \in Q$.)
Thus if $k$ is fixed,
by taking $n = k$, $l = p = 1$, and
$Q = \lbrace [y, y] \ : \ 0 \leq y \leq M-1 \rbrace$,
we can decide in polynomial time whether 
$\forall t, 0 \leq t \leq M-1, \ \exists x_1, x_2, \ldots, x_k \in \Zee$
such that $x_1, x_2, \ldots, x_k \geq 0$,
$\sum_{1 \leq i \leq k} x_i a_i = t$, and
$\sum_{1 \leq i \leq k} x_i \leq h$.
Now we use a binary search on $M$ to find the largest $M$ for which
the statement holds.  This gives
us $N_h (a_1, a_2, \ldots, a_k)$. 
Note that $Q$ is the intersection of four half-spaces.
Kannan's algorithm is quite complicated and this method is likely not
to be useful in practice.

\newcommand{\noopsort}[1]{} \newcommand{\singleletter}[1]{#1}

\end{document}